# Research note on a well posed integral used in Apery's proof for the irrationality of $\zeta(3)$

**Nikos Bagis**
Department of Informatics
Aristotele University of Thessaloniki Greece
bagkis@hotmail.com

**Abstract**
In this note we evaluate multiple integrals that play a crucial role in the theory of irrationality of zeta function.

By generalizing Beukers arguments in [1] we prove

**Proposition.**
Let $P_n(x) = \frac{1}{n!}\frac{d^n}{dx^n}\left(x^n(1-x)^n\right)$, and $v = 0,1,2,\ldots$, then

$$\int_0^1\int_0^1\cdots\int_0^1 \frac{(-\log(x_1x_2\cdot\ldots\cdot x_r))^v}{1-x_1x_2\cdot\ldots\cdot x_r}P_n(x_1)P_n(x_2)\cdot\ldots\cdot P(x_r)dx_1dx_2\ldots dx_r =$$

$$= \frac{d^v}{dz^v}\left(\sum_{k=0}^{\infty}\frac{1}{(z+n+k+1)^r}\prod_{j=0}^{n}\left(\frac{z-j+k+1}{z+j+k}\right)^r\right)_{z=0}$$

**Proof.** Set

$$c_v(n) = \int_0^1\int_0^1\cdots\int_0^1 \frac{(-\log(x_1x_2\cdot\ldots\cdot x_r))^v}{1-x_1x_2\cdot\ldots\cdot x_r}P_n(x_1)P_n(x_2)\cdot\ldots\cdot P_n(x_r)dx_1dx_2\ldots dx_r$$

then

$$\sum_{v=0}^{\infty}\frac{c_v(n)(-z)^v}{v!} =$$

$$= \sum_{v=0}^{\infty}\frac{(-z)^v}{v!}\int_0^1\int_0^1\cdots\int_0^1 \frac{(-\log(x_1x_2\cdot\ldots\cdot x_r))^v}{1-x_1x_2\cdot\ldots\cdot x_r}P_n(x_1)P_n(x_2)\cdot\ldots\cdot P_n(x_r)dx_1dx_2\ldots dx_r =$$

$$= \int_0^1\int_0^1\cdots\int_0^1\sum_{v=0}^{\infty}\frac{(-z)^v(-\log(x_1x_2\cdot\ldots\cdot x_r))^v}{v!(1-x_1x_2\cdot\ldots\cdot x_r)}P_n(x_1)P_n(x_2)\cdot\ldots\cdot P_n(x_r)dx_1dx_2\ldots dx_r =$$

$$= \int_0^1\int_0^1\cdots\int_0^1 e^{z\log(x_1x_2\cdot\ldots\cdot x_r)}\frac{1}{(1-x_1x_2\cdot\ldots\cdot x_r)}P_n(x_1)P_n(x_2)\cdot\ldots\cdot P_n(x_r)dx_1dx_2\ldots dx_r =$$

$$=\int_0^1\int_0^1\cdots\int_0^1 (x_1x_2\cdot\ldots\cdot x_r)^z \frac{1}{(1-x_1x_2\cdot\ldots\cdot x_r)} P_n(x_1)P_n(x_2)\cdot\ldots\cdot P_n(x_r)dx_1dx_2\ldots dx_r =$$

$$=\int_0^1\int_0^1\cdots\int_0^1 (x_1x_2\cdot\ldots\cdot x_r)^z \sum_{k=0}^{\infty}(x_1x_2\cdot\ldots\cdot x_r)^k P_n(x_1)P_n(x_2)\cdot\ldots\cdot P_n(x_r)dx_1dx_2\ldots dx_r =$$

$$=\sum_{k=0}^{\infty}\int_0^1\int_0^1\cdots\int_0^1 (x_1x_2\cdot\ldots\cdot x_r)^z (x_1x_2\cdot\ldots\cdot x_r)^k P_n(x_1)P_n(x_2)\cdot\ldots\cdot P_n(x_r)dx_1dx_2\ldots dx_r =$$

$$=\sum_{k=0}^{\infty}\left(\int_0^1 x^{z+k}P_n(x)dx\right)^r .$$

Hence we have $\displaystyle\sum_{v=0}^{\infty}\frac{c_v(n)(-z)^v}{v!}=\sum_{k=0}^{\infty}\left(\int_0^1 x^{z+k}P_n(x)dx\right)^r$ : (1)

by induction is

$$\int_0^1 x^{z+k}P_n(x)dx=\frac{1}{z+k+1}\prod_{j=0}^{n}\left(\frac{z-j+k+1}{z+k+j}\right) : (2)$$

From Taylor's expansion theorem we get the result. □

There is a interesting note here. If we want to use other polynomials except $P_n$ we set $R_n(x)=\sum_{k=0}^{n}a(n,k)x^k$ . Then again can be shown that:

$$\sum_{k=0}^{\infty}\left(\int_0^1 x^{z+k}R_n(x)dx\right)^r=\sum_{k=0}^{\infty}\left(\sum_{l=0}^{n}\frac{a(n,l)}{(l+k+z+1)}\right)^r$$

And a more general result read as

**Theorem.** Let $R_n(x)=\sum_{k=0}^{n}a(n,k)x^k$ be a family of polynomials for $v = 1,2,3\ldots$ we have

$$\int_0^1\int_0^1\cdots\int_0^1 \frac{(-\log(x_1x_2\cdot\ldots\cdot x_r))^v}{1-x_1x_2\cdot\ldots\cdot x_r}R_n(x_1)R_n(x_2)\cdot\ldots\cdot R_n(x_r)dx_1dx_2\ldots dx_r =$$

$$=\frac{d^v}{dz^v}\left(\sum_{k=0}^{\infty}\left(\sum_{l=0}^{n}\frac{a(n,l)}{(l+k+z+1)}\right)^r\right)_{z=0}$$

**Example.**
These generalizations with proper conditions give rise to some interesting results. For example if we consider the case $r = 3$, $v=2$ then we get that exists sequences $A(n)$, $B(n)$, $G(n)$, $d_n^5\in\mathbb{N}$ such that

$$\int_0^1\int_0^1\int_0^1 \frac{\log(x_1x_2x_3)^2}{1-x_1x_2x_3} P_n(x_1)P_n(x_2)P_n(x_3)dx_1dx_2dx_3 = \frac{A(n)\pi^4 + B(n)\zeta(5) + G(n)}{d_n^5}$$

with $d_n^{\,5} \mid lcm(1,2,3,...,n)^5$ and

$$\mathbb{N} \ni lcm(1,2,...,n)^5 \int_0^1\int_0^1\int_0^1 \frac{\log(x_1x_2x_3)^2}{1-x_1x_2x_3} P_n(x_1)P_n(x_2)P_n(x_3)dx_1dx_2dx_3 =$$

$$= lcm(1,2,...,n)^5 \frac{A(n)\pi^4 + B(n)\zeta(5) + G(n)}{d_n^5} << 1,$$

when $n >> 1$.

Hence if we assume that $\pi^4, \zeta(5)$ are both rational the quantity $Q$

$$Q = b_1b_2 lcm(1,2,...,n)^5 \frac{A(n)\pi^4 + B(n)\zeta(5) + G(n)}{d_n^5} =$$

$$= b_1b_2 lcm(1,2,...,n)^5 \frac{A(n)a_1/b_1 + B(n)a_2/b_2 + G(n)}{d_n^5} =$$

$$lcm(1,2,...,n)^5 \frac{d^2}{dz^2}\left(\sum_{k=0}^{\infty} \frac{1}{(z+n+k+1)^3} \prod_{j=0}^{n}\left(\frac{z-j+k+1}{z+j+k}\right)^3\right)_{z=0} < 1 \text{ for large } n$$

Hence $0 < Q(n) < 1$ and $Q(n) \in \mathbb{N}$, contradiction.

Thus for to prove that the numbers $\pi^4, \zeta(5)$, are not both rational we only have to prove that

$$e^{5n} \frac{d^2}{dz^2}\left(\sum_{k=0}^{\infty} \frac{1}{(z+n+k+1)^3} \prod_{j=0}^{n}\left(\frac{z-j+k+1}{z+j+k}\right)^3\right)_{z=0}$$ is as small as we want for large *n*.

**Note:** It is known that there exists constant *M* such that $lcm(1,2,...,n) < Me^n$.

In general for every *v*, *r* such that

$$e^{(r+v)n} \frac{d^v}{dz^v}\left(\sum_{k=0}^{\infty} \frac{1}{(z+n+k+1)^r} \prod_{j=0}^{n}\left(\frac{z-j+k+1}{z+j+k}\right)^r\right)_{z=0}$$ is as small as we want, for large *n*,

then the constants involving in

$$\int_0^1\int_0^1\cdots\int_0^1 \frac{(-\log(x_1x_2\cdot...\cdot x_r))^v}{1-x_1x_2\cdot...\cdot x_r} P_n(x_1)P_n(x_2)\cdot...\cdot P(x_r)dx_1dx_2...dx_r$$

must have among them at least one irrational.

In general the same things hold with an arbitrary polynomial *R*, with integer coefficients.